\documentclass[11pt]{article}
\usepackage{mathrsfs}
\usepackage{bbm}

\usepackage{verbatim}

\usepackage{amsthm}
\usepackage{amsmath}
\usepackage{amssymb}

\usepackage{graphicx}
\usepackage{epstopdf}
\usepackage{titling}

\usepackage{enumerate}

\usepackage{algpseudocode}
\usepackage{algorithm}
\usepackage{algorithmicx}
\usepackage{tikz}

\newenvironment{pf}
      {\medskip\noindent{\bf Proof.}\hspace{1mm}}
      {\hfill$\Box$\medskip}

\newtheorem{theorem}{Theorem}[section]

\newtheorem{lemma}[theorem]{Lemma}

\newtheorem{conjecture}[theorem]{Conjecture}

\begin{document}

\title{On the maximum induced density of directed stars\\ and related problems}

\author{Hao Huang\thanks{Institute for Advanced Study, Princeton, NJ 08540. Email: {\tt huanghao@math.ias.edu}. Research supported in part by NSF grant DMS-1128155.}}
\date{}
\maketitle
\setcounter{page}{1}
\vspace{-2em}
\begin{abstract}
Let $k \ge 3$ be an integer, we prove that the maximum induced density of the $k$-vertex directed star in a directed graph is attained by an iterated blow-up construction. This confirms a conjecture by Falgas-Ravry and Vaughan, who proved this for $k=3, 4$. This question provides the first known instance of density problem for which one can prove extremality of an iterated blow-up construction. We also study the inducibility of complete bipartite digraphs and discuss other related problems.
\end{abstract}

\section{Introduction}
In modern extremal combinatorics, a substantial number of problems study the asymptotic relations between densities of subgraphs, and can be formulated in the following language.
Given a family $\mathcal{F}$ of graphs and another graph $H$, define the \textit{Tur\'an $H$-number} of $\mathcal{F}$ to be the maximum number of induced copies of $H$ in a $\mathcal{F}$-free graph on $n$ vertices, and denote it by $\textrm{ex}_H(n, \mathcal{F})$.
We also denote by $\pi_{H}(\mathcal{F})$ the limit of the maximum induced density of $H$ in a $\mathcal{F}$-free graph when the number of vertices tends to infinity. Similar definitions can be as well made in the setting of $r$-uniform hypergraph, directed graph, and so forth.
When $H$ is a single edge, $\pi_{e}(\mathcal{F})$ is just the classical Tur\'an density. It has been a long-standing open problem in extremal combinatorics
to understand these densities for families of hypergraphs and directed graphs.
For results and techniques, we refer the readers to the survey \cite{keevash_survey}.

On the other hand, when $\mathcal{F}=\emptyset$, $\pi_H(\emptyset)$ studies the maximum induced density of $H$ in arbitrary graph, and is known
as the inducibility of $H$. Although there are various works \cite{behj_95, bnt_86, brown_sidorenko, exo_86, hirst_11, pg_75} on the inducibility of graphs, there are relatively fewer results for directed graphs. Sperfeld \cite{sperfeld} studied the inducibility of some digraphs on three vertices.
Falgas-Ravry and Vaughan \cite{turan_3-graph} determined $\pi_{\vec{S}_3}(\emptyset)$ and $\pi_{\vec{S}_4}(\emptyset)$ by flag algebra. Here
$\vec{S}_k$ is the directed star on $k$ vertices, with one vertex being the center, and $k-1$ edges oriented away from it.
They further made the following conjecture that the extremal digraph having the maximum induced density of $\vec{S}_k$ is always an unbalanced blow-up of $\vec{S}_2$ iterated inside one part.

\begin{conjecture}\label{conj_sk}
For every integer $k \ge 3$, $$\pi_{\vec{S}_k}(\emptyset) = \alpha_k = \max_{0 \le x \le 1} \frac{kx(1-x)^{k-1}}{1-x^k}.$$
Assume the maximum is attained by $x=x_k$, then the extremal configuration is constructed by starting with two parts $|A|=x_k n$ and $|B|=(1-x_k)n$, adding all the edges oriented from $A$ to $B$, and iterating this process inside $A$.
\end{conjecture}

As we mentioned earlier, the proofs for cases $k=3$ and $k=4$ employ the method of flag algebra developed by Razborov \cite{razborov} and are partly computer assisted. However since the search space and running time grow exponentially in $k$, a different approach may be needed for large $k$.
It is also worth mentioning that for most Tur\'an-type problems studying densities of subgraphs, if the conjectured extremal graph comes from such iterated construction instead of a simple blow-up of small graph, usually we do not know how to obtain an exact bound. For example, the Tur\'an density of $K_4^{-}$ (the unique $3$-graph on $4$ vertices with $3$ edges) was conjectured to be achieved by the iterated blow-up of certain $6$-vertex $3$-graph by Frankl and Fur\"edi \cite{ff_84} and is still open.
Another example is the special case of the well-known Caccetta-H\"aggkvist conjecture \cite{ch_78}: every $n$-vertex digraph with minimum outdegree at least $n/3$ contains a directed triangle. Its difficulty probably lies in the fact that the iterated blow-up of a directed $4$-cycle is one of the conjectured extremal examples. To the best of our knowledge, the case $k=3$ and $4$ of Conjecture \ref{conj_sk} are probably the only examples that the exact bound has been proved for an iterated construction, which leads us to believe there exists a simpler and more human-readable proof. Actually we are able to apply certain operations on digraphs, and reduce it to an optimization problem, and verify this conjecture for every directed star $\vec{S}_k$ for $k \ge 3$.

The rest of this short paper is organized as follows. In Section \ref{section_mainproof} we give a complete proof of Conjecture \ref{conj_sk}. Section \ref{section_gen} discusses the inducibility for complete bipartite digraphs. In the concluding remarks, we mention some related problems and possible future directions for research.

\section{Main proof} \label{section_mainproof}
In this section we will give a proof of Conjecture \ref{conj_sk}. Our proof is inspired by that of \cite{brown_sidorenko}.
Assume $D_n$ is the extremal directed graph on $n$ vertices which has the maximum number of induced copies of of $\vec{S}_k$.
We define an equivalence relation on its vertex set $V(D_n)$ as follows: $u \sim v$ iff they have the same in- and out-neighbordhoods, i.e. $N^+(u)=N^+(v)$ and $N^-(u)=N^-(v)$. This equivalence relation naturally
partitions the vertices of $D_n$ into the following equivalence classes: $V(D_n)=V_1 \cup \cdots \cup V_m$, where each $V_i$ induces an empty digraph.
From the definition, between two classes $V_i$ and $V_j$ there are three possible scenarios:
$(i)$ all the edges are oriented from $V_i$ to $V_j$; $(ii)$ all the edges are oriented from $V_j$ to $V_i$; $(iii)$ there is no edge between $V_i$ and $V_j$.
We claim that a sequence of operations can be applied on $D_n$ such that
the number of induced copies of $\vec{S}_k$ does not decrease, and in the resulted digraph case $(iii)$ never occurs.

\begin{lemma}\label{lemma_non_empty}
Given a directed graph $D_n$ with equivalence classes $V_1, \cdots, V_m$, and $1 \le i \neq j \le m$. If there is no edge between $V_i$ and $V_j$, we can merge $V_i$ and $V_j$ into one equivalence class without decreasing the number of induced copies of $\vec{S}_k$.
\end{lemma}
\begin{pf}
Assume $|V_i|=x$ and $|V_j|=y$. Denote by $D'_n$ the digraph formed by moving vertices between $V_i$ and $V_j$ and changing their neighborhoods accordingly, with $|V'_i|=z$, $|V'_j|=x+y-z$. Let $N_{00}$ be the number of induced copies of $\vec{S}_k$ in $D$ not involving vertices in $V_i$ or $V_j$; $N_{10}$ be the number of induced copies of $\vec{S}_k$ using vertices from $V_i$ but not vertices from $V_j$; $N_{01}$ be the number of induced copies of $\vec{S}_k$ using vertices from $V_j$ but not vertices from $V_i$; and finally $N_{11}$ be the number of induced copies of $\vec{S}_k$ using vertices from both $V_i$ and $V_j$. Obviously the total number of induced copies of $\vec{S}_k$ is equal to $N_{00}+N_{10}+N_{11}+N_{11}$. Similarly we also define the parameters $N'_{00}, N'_{10}, N'_{01}, N'_{11}$ for $D'_n$.

Note that $N'_{00}=N_{00}$ since only the adjacencies involving vertices in $V_i \cup V_j$ might be changed. We also have $N_{11}=N'_{11}$.
Consider an induced copy of $\vec{S}_k$ containing $v_i \in V_i$, and $v_j \in V_j$. Since there is no edge between the two parts $V_i$ and $V_j$, $v_i$ is not
adjacent to $v_j$. Therefore both of them are the leaves of $\vec{S}_k$. Because all the leaves are equivalent in $\vec{S}_k$,
moving vertices between $V_i$ and $V_j$ does not change the value of $N_{11}$.

Next we show that $z$ can be chosen such that $N'_{10}+N'_{01} \ge N_{10}+N_{01}$.
Denote by $s_l$ the number of $(k-l)$-vertex sets $S$ in $[n]\backslash (V_i \cup V_j)$ such that $S$ together with any $l$ vertices in $V_i$ induce a copy of $\vec{S}_k$. Similarly let $t_l$ be the number of $(k-l)$-vertex sets $T$ such that $T$ together with any $l$ vertices in $V_j$ induce a copy of $\vec{S}_k$. Then by the definition of equivalence class, we have
\begin{align*}
N_{10}=\sum_{l=1}^{k} \binom{x}{l} s_l, &~~~~~~~ N_{01}=\sum_{l=1}^{k} \binom{y}{l} t_l,\\
N'_{01}=\sum_{l=1}^{k} \binom{z}{l} s_l, &~~~~~~~ N'_{01}=\sum_{l=1}^{k} \binom{x+y-z}{l} t_l.
\end{align*}

It is not difficult to verify that $\binom{z}{l}$ and $\binom{x+y-z}{l}$ are both convex functions in the variable $z$. Therefore one could merge these two equivalence classes $V_i$ and $V_j$ by taking either $z=0$ or $z=x+y$ in the new digraph $D'_n$, such that $N'_{01}+N'_{10} \ge N_{01}+N_{10}$.
\end{pf}
\medskip

Note that after merging vertices in Lemma \ref{lemma_non_empty}, the number of equivalence classes decreases, so this process stops after a finite number of steps. We may assume that in the extremal digraph $D_n$ with equivalent classes $\{V_1, \cdots, V_m\}$, and any $i \neq j$, either there is a complete bipartite digraph with every edge oriented from $V_i$ to $V_j$, or from $V_j$ to $V_i$. We denote them by $V_i \rightarrow V_j$ and $V_j \rightarrow V_i$ respectively. Assume $|V_i|=w_i n$ where $\sum_{i=1}^m w_i=1$. The induced density of $\vec{S}_k$ in $D_n$ is equal to
$$\dfrac{1}{\binom{n}{k}}\sum_{V_i \rightarrow V_j} w_i n \binom{w_j n}{k-1} = \sum_{V_i \rightarrow V_j} k w_i w_j^{k-1} + o(1).$$
Since $\pi_{\vec{S}_k}(\emptyset)$ is the limit of the maximum induced density when $n$ tends to infinity, we can neglect the $o(1)$ term here.
Without loss of generality, we may assume that $w_1 \ge w_2 \ge \cdots \ge w_m$ by reordering $\{V_i\}$. If $V_i \rightarrow V_j$ in $D_n$ for some $i<j$, then by changing the orientation of this complete bipartite digraph, the induced density increases
by $w_j w_i^{k-1} -w_i w_j^{k-1} \ge 0$. Therefore we can assume $V_j \rightarrow V_i$ for any $i<j$. Basically speaking we obtain $D_n$ to be the unbalanced blow-up of a transitive tournament, and the induced density of $\vec{S}_k$ in $D_n$ is now equal to
\begin{align*}
f_m(w_1, \cdots, w_m) &= k \sum_{1 \le i <j \le m} w_i^{k-1} w_j\\
&=k \cdot \left(w_1^{k-1}(w_2+\cdots+w_{m})+w_{2}^{k-1}(w_3+ \cdots+ w_{m}) + \cdots +w_{m-1}^{k-1} w_m\right)
\end{align*}
Let $F_m(x) = \max f_m(w_1, \cdots, w_m)$ subject to $\sum_i w_i=x$ and $w_i \ge 0$, then $\pi_{\vec{S}_k}(\emptyset) = \limsup_{m \rightarrow \infty} F_m(1)$.
Because $f_m$ is a homogeneous polynomial of degree $k$, we have $F_m(x)=F_m(1)x^k$ and thus
\begin{align} \label{homo}
F_m(1) &= \max_{0 \le w_1 \le 1} kw_1^{k-1}(1-w_1) + F_{m-1}(1-w_1) \nonumber\\
 &= \max_{0 \le w_1 \le 1} kw_1^{k-1}(1-w_1) + (1-w_1)^k F_{m-1}(1).
\end{align}
Taking $w_1=0$ in \eqref{homo} shows that $F_m(1) \ge F_{m-1}(1)$. Due to the fact that the induced density can never be greater than $1$, $\{F_m(1)\}$ is a bounded monotone non-decreasing sequence and thus converges to a limit, denoted by $\alpha_k$. Let $m \rightarrow \infty$ in \eqref{homo}, we have
$$\alpha_k = \max_{0 \le x \le 1} kx^{k-1}(1-x) + (1-x)^k \alpha_k$$
Let $\beta_k = \max_{0 \le x \le 1} \dfrac{kx^{k-1}(1-x)}{1-(1-x)^k}$, we now prove that $\alpha_k=\beta_k$. Since $\dfrac{kx^{k-1}(1-x)}{1-(1-x)^k}$ is continuous and bounded on the compact set $[0,1]$, $\beta_k = \dfrac{ky^{k-1}(1-y)}{1-(1-y)^k}$ for some $y \in [0,1]$ and thus
$$\alpha_k \ge ky^{k-1}(1-y) + (1-y)^k \alpha_k = (1-(1-y)^k) \beta_k +(1-y)^k \alpha_k,$$
which implies that $\alpha_k \ge \beta_k$. On the other hand, suppose $z \in [0,1]$ maximizes $kx^{k-1}(1-x) + (1-x)^k \alpha_k$, then
$\alpha_k = kz^{k-1}(1-z)+(1-z)^k \alpha_k,$ and
$$\beta_k \ge \dfrac{kz^{k-1}(1-z)}{1-(1-z)^k}=\alpha_k.$$
Therefore $$\pi_{\vec{S}_k}(\emptyset) = \alpha_k= \beta_k = \max_{0 \le x \le 1} \dfrac{kx^{k-1}(1-x)}{1-(1-x)^k} = \max_{0 \le x \le 1} \dfrac{kx(1-x)^{k-1}}{1-x^k}.$$
Suppose $x_k\in[0, 1]$ maximizes $kx(1-x)^{k-1}/(1-x^k)$, from the above proof, we can see that the bound is obtained uniquely by the infinite sequence
$w_i = x^{i-1}(1-x_k)$, $i=1, 2, \cdots$, which corresponds to the iterated blow-up of $\vec{S}_2$ and concludes the proof of Conjecture \ref{conj_sk}.

\section{Inducibility of complete bipartite digraph} \label{section_gen}
Note that the proof of Lemma \ref{lemma_non_empty} also works for any digraph $H$ in which every two vertices are in the same equivalence class whenever
they are not adjacent. In particular, it works for the problem of maximizing the induced density of complete bipartite digraph $\vec{K}_{s,t}$. Here $V(\vec{K}_{s,t})=[s+t]$, and the edge set consists of edges from $i$ to $s+j$ for every $1 \le i \le s, 1 \le j \le t$.
When $s=1$, this corresponds to the directed star $\vec{S}_{t+1}$.
Since $\pi_{\vec{K}_{s,t}}(\emptyset)=\pi_{\vec{K}_{t,s}}(\emptyset)$ by flipping the orientation of every edge, we can assume $s \leq t$. Similarly as before,
if the equivalence classes $V_1, \cdots, V_m$ in the extremal digraph $D_n$ has sizes $|V_i|=w_i n$ with $w_1 \ge w_2 \ge \cdots \ge w_m$, then the maximum induced density $\pi_{\vec{K}_{s,t}}(\emptyset)$ is equal to $\lim_{m \rightarrow \infty} F_m$, where
$F_m$ is the maximum of $f_m(w_1, \cdots, w_m)$ subject to $\sum_{i=1}^m w_i=1$ and $w_i \ge 0$, with
$$f_m(w_1, \cdots, w_m) = \binom{s+t}{s} \sum_{1 \le i<j \le m} w_i^t w_j^s.$$

\begin{theorem}\label{max_kst}
For integers $2 \le s \le t$,
$$\pi_{\vec{K}_{s,t}}(\emptyset) = \binom{s+t}{s} \left(\dfrac{s}{s+t}\right)^s \left(\dfrac{t}{s+t} \right)^t,$$
which is achieved by the balanced blow-up $\vec{K}_{\frac{s}{s+t}n, \frac{t}{s+t}n}$ of $\vec{K}_{s,t}$, when $n \rightarrow \infty$.
\end{theorem}
\begin{pf}
First by taking $m=2$ and $w_1 = \frac{t}{s+t}$ and $w_2 = \frac{s}{s+t}$, we have
$$\pi_{\vec{K}_{s,t}}(\emptyset) = \lim_{m \rightarrow \infty} F_m \ge F_2 \ge \binom{s+t}{s} \left(\dfrac{s}{s+t}\right)^s \left(\dfrac{t}{s+t} \right)^t.$$
On the other hand,
\begin{align} \label{ineq_kst}
f_m(w_1, \cdots, w_m) &= \binom{s+t}{s} \sum_{1 \le i<j \le m}w_i^t w_j^s \nonumber\\
&= \binom{s+t}{s} \left(w_1^t\left(\sum_{i=2}^m w_i^s\right)+\sum_{2 \le i<j \le m} w_i^t w_j^s\right) \nonumber\\
&\le \binom{s+t}{s} \left(w_1^t\left(\sum_{i=2}^m w_i^s\right)+\sum_{2 \le i<j \le m} w_1^{t} w_i w_j^{s-1}\right)\\
&\le \binom{s+t}{s} w_1^t (w_2+ \cdots + w_m)^s.\nonumber
\end{align}
The first inequality is because $w_i \le w_1$ and $w_j \le w_1$. The second inequality follows from the fact that the coefficient of $w_i w_j^{s-1}$ in the expansion of $(w_2+ \cdots+ w_m)^s$ is equal to $s$, which is greater than the coefficient $1$ of the corresponding term in the left hand side, whenever $s \ge 2$. It follows from inequality \eqref{ineq_kst} that $F_m \le F_2$. By elementary calculus, one can easily show that $w_1^{t}w_2^{s}$ is maximized when $w_1=\frac{t}{s+t}$ and $w_2=\frac{s}{s+t}$, which finishes the proof of Theorem
\ref{max_kst}.
\end{pf}

\section{Concluding remarks}
\begin{list}{\labelitemi}{\leftmargin=1em}
\item In \cite{turan_3-graph}, the authors mention that for any given digraph $D$, an auxiliary $3$-uniform hypergraph $G(D)$ can be defined by setting
$xyz$ to be a $3$-edge whenever $\{x, y, z\}$ induces a copy of $\vec{S}_3$ in $D$. It is not hard to check that $G(D)$ is always a $C_5$-free $3$-graph. Here $C_5$ refers to the tight cycle on $5$ vertices, whose edges are $(123)$, $(234)$, $(345)$, $(451)$, and $(512)$. Mubayi and R\"odl conjectured that the Tur\'an number $\pi(C_5)$ is equal to $2\sqrt{3}-3$, with exactly the same iterated construction in the $\vec{S}_3$ problem.
The result $\pi_{\vec{S}_3}(\emptyset)=\max_{0 \le x \le 1} 3x(1-x)^2/(1-x^3) = 2\sqrt{3}-3$ settles the special case when the $3$-graph has the form $G(D)$ from a digraph $D$.

\item Sperfeld \cite{sperfeld} studies the maximum induced density of some small digraphs, and in particular he proved that $\pi_{\vec{C}_3}(\emptyset) = 1/4$ with the extremal example including the random tournament and the iterated blow-up of $\vec{C}_3$, and conjecture that $\pi_{\vec{C}_4}(\emptyset)$ is achieved by the iterated blow-up of $\vec{C}_4$. It would be of great interest to develop new techniques to attack this problem, since the solution of this problem might as well provide insights into solving the Caccetta-H\"aggkvist conjecture.

\item Although obtaining a general solution to the graph or digraph inducibility problem seems to be difficult,
Hatami, Hirst and Norine  \cite{hatami_hirst_norine} showed that for a given graph $H$, the $n$-vertex graph $G$ containing the most number of induced copies of sufficiently large balanced blow-up of $H$, is itself essentially a blow-up of $H$. It would be interesting if similar results can be proved for the inducibility of directed graphs, which may also involve some iterated blow-ups.
\end{list}

\vspace{0.25cm}
\noindent \textbf{Acknowledgement.} The author would like to thank Benny Sudakov for his valuable comments and discussions.

\end{document}